\input amstex
\documentstyle{amsppt}
\input bull-ppt
\keyedby{bull507e/kmt}
\ratitle
\topmatter
\cvol{31}
\cvolyear{1994}
\cmonth{July}
\cyear{1994}
\cvolno{1}
\cpgs{68-74}
\title Open sets of diffeomorphisms having two attractors, 
\\
each with an everywhere dense basin \endtitle
\shorttitle 
{Two Attractors, Each With an Everywhere Dense Basin}
\author Ittai Kan \endauthor
\address Department of Mathematics, George Mason 
University, Fairfax,
Virginia 22030\endaddress
\ml ikan\@mason1.gmu.edu \endml 
\date March 23, 1993\enddate
\abstract
We announce the discovery of a diffeomorphism of a 
three-dimensional manifold with boundary which has 
two disjoint attractors. Each attractor attracts a set of 
positive
$3$-dimensional Lebesgue measure whose points of Lebesgue 
density are
dense in the whole manifold. 
This situation is stable under small perturbations.
\endabstract 

\thanks Supported in part by AFOSR grant 
F49620-92-J-0033\endthanks
\endtopmatter

\document

\heading 
1. Introduction
\endheading 

Through the numerical experiments and theoretical works of 
many researchers 
we have gained some understanding of a dynamical system  
which is  
``typical'' in the ``standard'' setting; that is, the 
dynamical 
system $f:M \to M$ is a smooth map 
and $M$ is a smooth $n$-dimensional manifold. 
Special circumstances change the meaning of ``typical'' 
and add 
possibilities to the list of what we are likely to see, 
say, by
requiring $f$ to be in some particular subset of smooth 
maps on $M$. 
For example, 
if $f$ is symplectic, then we see the KAM tori and other 
phenomena of 
Hamiltonian dynamics; or if $f$ has symmetries, 
we are likely to see ``symmetry-breaking'' 
bifurcations.  The setting we examine here is 
where $f$  
maps
some fixed submanifold $K \subset M$ of dimension
$k < m$ to itself ($f$ is {\it invariant} on $K$).
That is, we discuss some new phenomena which are
``typical'' in the sense that they occur for some large 
set of 
perturbations of $f$ which continue to map $K$ to itself. 
This is the 
natural setting for modeling  
physical systems with constraints or 
reflectional symmetries.

We find some unexpected new limiting behaviors in the new 
setting. 
More precisely, 
the $\omega$-{\it limit} set of a point ${\bold x} \in M$ 
is the 
set $\omega({\bold x}) = \bigcap_{k} \overline{\bigcup_{l 
\ge k} f^{l}({\bold x})}$, the forward limit set
of ${\bold x}$ under iteration by $f$. We say that a 
compact set $A$ is 
an {\it attractor}  if 
$\{ {\bold x} \in M : \omega({\bold x}) = A \}$ has  
positive $n$-dimensional Lebesgue
measure and 
the {\it basin of attraction} of $A$ is ${\Cal B}(A) = \{ 
{\bold x} \in M : \omega({\bold x}) \subset A \}$.

In the standard 
setting the 
basins of  
attraction typically seem to be open sets, although the 
boundaries of these 
basins are often fractal sets \cite{1}. Certainly, in the 
standard setting, 
it has long been known that 
the only structurally stable attractors are the Axiom A 
attractors, and these 
have open basins of attraction; see, for example, \cite{2}.
In the setting where we have a fixed invariant manifold 
$K$, we find
that attractors may have basins with some radically new 
features.
We find that some maps have several coexisting attractors 
whose basins are
measure-theoretically dense in each other. 
That is, 
any open set which intersects one basin in a set of 
positive 
measure also intersects each of the other basins in a set 
of positive 
measure. We call such basins {\it intermingled}.

In \cite{3}   
numerical evidence for intermingled basins is presented in 
the case of 
a particular quadratic map of the plane. 
The main result of this paper is that there are maps for 
which 
intermingling can be rigorously shown to occur and, in 
fact, to persist 
under perturbation. 

\proclaim{Theorem 1} For every $k \geq 1 + \alpha$, 
there is an open set of $C^k$
diffeomorphisms of the thickened torus with boundary
$(T^2 \times I)$ for which there are
two coexisting attractors whose basins are intermingled and
the union of whose basins has full measure.  
\endproclaim
		   
In fact, with choice of coordinates $T^2 =
S^1 \times S^1 = [0,\,1) \times [0,\,1)$ and
$I = [0,\,1]$ we find that a small $C^k$ 
neighborhood of the $C^{\infty}$ 
diffeomorphism  
$$f(x,\,y,\,z) = (3x+y,\,2x+y,\,z+\cos(2\pi 
x)({z\over32})(1-z)
) \tag 1$$
satisfies the theorem. 
That is, for any diffeomorphism in this neighborhood, 
the two boundary tori, $A = T^2 \times \{ 0 \}$ and 
$B = T^2 \times \{1 \}$, are each attractors; the union of 
their 
basins has full $3$-dimensional Lebesgue measure; and 
their basins are 
intermingled.

A result similar to Theorem 1 can be obtained in two 
dimensions if we 
abandon the requirement of invertibility.  For the annulus 
we have the 
following theorem. 

\proclaim{Theorem 2} For every $k \geq 1+ \alpha$, 
there is an open set of $C^k$
maps of the annulus $(S^1 \times I)$ 
which send the boundary to the boundary
for which there are
two coexisting attractors whose basins are intermingled and
the union of whose basins has full measure.  
\endproclaim

With choice of coordinates $S^1 = [0,\,1) $ and
$I = [0,\,1]$, we find that
a small $C^k$ neighborhood of the $C^{\infty}$ 
map   
$$f(x,\,z) =(3x,\,z+\cos(2\pi x)({z\over32})(1-z)
)\tag 2$$
satisfies the theorem. 
That is, for any map in this neighborhood, 
the two boundary circles, $A = S^1 \times \{ 0 \}$ and 
$B = S^1 \times \{1 \}$, are each attractors; the union of 
their 
basins has full $2$-dimensional Lebesgue measure; and 
their basins are 
intermingled.

\rem
{Remark {\rm1}} Although the basins for these attractors 
are 
measure-theoretically large, they are small topologically
(they are meager, the countable union of nowhere dense 
sets). 
It is easy to show that this is generally true
in the case of a continuous map where there are two
disjoint attractors whose basins are
intermingled.
\endrem 

\rem
{Remark {\rm2}} There are many other definitions of 
attractor and basin in use 
today \cite{4, 5}. The theorems of this paper 
hold for most of the measure-theoretic
definitions. 
\endrem 

The proofs of Theorems 1 and 2 are long and technical. In 
Section 2  
we give a proof of Theorem 2 restricted to the case where 
$f$ is the map given by (2). This proof is elementary and 
provides 
insight into the general proofs.
In Section 3 we discuss the 
modifications needed to prove the general case of Theorems 
1 and\ 2.

\heading 
2. An example on the annulus
\endheading 

The proof of Theorem 1 is 
similar to the proof of Theorem 2 but contains additional 
technicalities which may obscure the main ideas of this 
paper.  
In addition, the $2$-dimensional dynamics of  
the map of the annulus 
given by (2) are much easier to visualize than the 
the map of the thickened torus 
given by (1); so in this section we show that for the map 
of the annulus 
given by (2), the two boundary circles are each attractors, 
the union of their basins has full $2$-dimensional 
Lebesgue measure, 
and their basins are intermingled. 

Our coordinates for the annulus will be 
$(x,\, z) \in [0,1) \times [0,1]$ with $(0,\,z)$ 
identified with $(1,\,z)$. 
We note that for $f$, the $x$-coordinate function $f_x$ 
depends only on
the $x$-coordinate, and thus $f$ maps vertical line 
segments to 
vertical line segments.  
We use $m$ to denote $2$-dimensional Lebesgue measure and  
$m_x$ to denote the measure on curves which is locally 
$1$-dimensional 
Lebesgue measure of the orthogonal projection to the 
$x$-axis.

We note that the points $(0,\,0),$ $(0,\,1),$ $({1\over 
2},\,0),$ 
and $({1\over 2},\,1)$ are each fixed points of ${f}.$  
Also, the Jacobian matrix has the lower triangular form 
$$\split 
J{f}(x,\,z) 
&=\pmatrix [D_x{f_{x}}](x,\,z)
& [D_z{f_{x}}](x,\,z) \\
[D_x { f_{z}}](x,\,z)
&[D_z {f}_{z}](x,\,z)\endpmatrix \\ 
&= \pmatrix 3 & 0 \\
-{\pi z \over 16}(1-z) \sin(2 \pi x) &
1+{(1-2z) \over 32}\cos(2 \pi x)\endpmatrix . 
\endsplit \tag 3$$ 

\subheading{Lyapunov exponents and stable manifolds} 
We now review some basic facts and definitions from 
the theory of dynamical systems. In particular, we discuss 
Lyapunov exponents and stable manifolds.

The Lyapunov exponent for a vector ${\overrightarrow v}$
with base point $(x,\,z)$ is the exponential rate 
at which the vector length grows under iteration by $f$. 
That is,
$${\lambda}({\overrightarrow v}) =
\lim_{n \to \infty} {1\over n} \ln{\| J f^n (x,\,z) 
{\overrightarrow v}\|}$$
if the limit exists and where $f^n$ denotes the $n$-th 
iterate of $f$. 
Each point in ${\bold x} \in A$
has at most two Lyapunov 
exponents. 
One is called the {\it normal Lyapunov 
exponent} and measures the exponential rate of expansion 
in the direction 
perpendicular to the boundary circle $A$. By the diagonal 
form of the 
Jacobian on $A$ this exponent is   
$$\lambda_{\perp}({\bold x})=
\lim_{n \to \infty}{1\over n}\sum^{n}_{k=1}{
\ln{\vert [D_{z}f_{z}](f^k({\bold x}))\vert}},$$  
if the limit exists.

In order to organize which points are attracted to which 
attractor, 
the following definitions are  
useful. 
The {\it stable manifold} of a point ${\bold x}$ is the 
set of points which converge 
to the orbit of ${\bold x}$, that is, $W^{\roman s}({\bold 
x}) = 
\{ {\bold y} \in M  : \lim_{n \to \infty} \| 
f^{n}({\bold x}) - f^{n}({\bold y}) \| = 0\}.$ 
Since ${f}$ maps vertical segments to vertical segments and 
the restriction of $f$ to a single vertical segment is 
invertible,
it is easy to see that if the point $({\hat x},\,{\rho})$ 
is in 
the stable manifold for a point ${\bold x} =( x,\,0) \in 
A$, 
then every point
$({\hat x},\,z)$ with
$z \le {\rho}$ is also in the stable manifold of 
the point ${\bold x}$. Thus the stable manifold\ of 
the point $\bold x$ is 
a union of (possibly trivial)
vertical line segments whose bottoms terminate in $A$. 
The component of
$W^{\roman s}({\bold x})$ which contains ${\bold x}$ is
called the {\it immediate stable
manifold} of ${\bold x}$ and is denoted by 
$W^{\roman{is}}({\bold x})$. 
Since $f$ is expanding in the $x$-direction
and is $3$-to-$1$, we see that the stable manifold of 
${\bold x}$
is a union of inverse images of forward images of the 
immediate
stable manifold of ${\bold x}$. That is, $W^{\roman 
s}({\bold x}) =
\bigcup_{n \ge 0}
{{f}^{-n}({f}^{n}(W^{\roman{is}}({\bold x})))}$. 

\subheading{Lemmas}
We now show that $A$ is on the average contracting in the 
normal direction. 
This will  
imply that $A$ is an attractor and that the immediate 
stable manifold 
of almost every point in $A$ is a vertical line segment of 
positive 
length. 
The analogous results for the other boundary circle $B$ 
follow by 
symmetry.  

\proclaim{Lemma {\rm2.1}} For $m_x$-almost every  
point ${\bold x} \in A$ the normal Lyapunov exponent
is negative. 
\endproclaim

\demo{Proof of Lemma {\rm2.1}} Since $1$-dimensional 
Lebesgue 
measure $m_x$ on $A$ is ergodic and 
invariant under $f$, the ergodic theorem tells us that  
almost every ${\bold x}= (x,\,0) \in A$
has a normal Lyapunov exponent equal to 
$$\lambda_{\perp} = \int^{1}_{0} \ln \vert 
[D_{z}{f_z}](x,\,0) \vert dx = 
\int^{1}_{0} \ln \vert 1 +
{\cos(2 \pi x)\over 32} \vert\, dx .$$
By symmetry we have
$$\align
\lambda_{\perp} &= 2\int^{1 \over 4}_{0}  
\ln \vert 1 -
{\cos(2 \pi x)\over 32} \vert + \ln \vert 1 +
{\cos(2 \pi x)\over 32} \vert\, dx\\
&= 
2\int^{1 \over 4}_{0}  
\ln \vert 1 -
\bigl({\cos(2 \pi x)\over 32} \bigr)^2 \vert\, dx < 
0.\endalign
$$
\enddemo

We now show that for every point in $A$ whose normal 
Lyapunov exponent is 
negative the immediate stable manifold has positive length. 
Since $m_x$-almost every point of $A$ has a normal 
Lyapunov exponent 
equal to $\lambda_{\perp} < 0$ and all immediate stable 
manifolds are 
vertical line segments, we then easily see that $m({\Cal 
B}(A)) >0$.

\proclaim{Lemma {\rm2.2}} If 
${\bold {\hat x}} = ({\hat x},\,0) \in A$ and 
$\lambda_{\perp}
({\bold {\hat x}})= \lambda_{\perp}<0$, 
then $W^{\roman{is}}({\bold {\hat x}})$ has positive length.
\endproclaim

\proclaim{Corollary} The basin of $A$ has positive 
$2$-dimensional Lebesgue 
measure.
\endproclaim

\demo{Proof of Lemma {\rm2.2}} Assume ${\bold {\hat x}}$ 
is as in the hypothesis of the lemma. 
Choose a positive $c$ sufficiently small so that for all 
$x$ we have that 
$z<c$ implies 
$\vert \ln{\vert [D_zf_z](x,\,0)\vert} - \ln{\vert 
[D_zf_z](x,\,z)\vert} \vert 
< \vert {\lambda_{\perp} \over 4} \vert$. 
Since $\lambda_{\perp}({\bold {\hat x}}) < 0$, the 
quantity  $\tau({\hat x})$ defined as follows
$$\tau({\hat x}) = \sup_{n \geq 0}\{ {-3n\lambda_{\perp} 
\over 4}+
\sum_{k=0}^{n}  
{\ln \vert [D_{z}{f_z}]f^{k}({\bold {\hat x}}) \vert } \}, 
$$ 
is finite. 

By the Mean Value Theorem and the Chain Rule for 
differentiation, 
we see that if $f^{k}_{z}({\hat x},\,z) < c$ for all 
$0 \leq k \leq n$, then 
$$\ln{\vert f^{n+1}_{z}({\hat x},\,z) \vert } < 
\ln{(z)} - \vert {n\lambda_{\perp} \over 4} \vert+
{\sum_{k=0}^{n}  
{\ln \vert [D_{z}{f_z}]f^{k}({\bold {\hat x}}) \vert}} \leq 
\ln{(z)} + \tau({\bold {\hat x}}) +{n \lambda_{\perp} 
\over 2}. $$ 
In fact, if $z$ is chosen to be less than ${c \exp{( 
-\tau({\hat x}))}},$ 
then $f_z^{n+1}({\hat x},\,z) < c \exp{({n \lambda_{\perp} 
\over 2})}<c$. 
Thus by induction we see that the length of $ 
W^{\roman{is}}({\bold {\hat x}}) $ 
is at least ${c \exp{( -\tau({\hat x}))}}$. 
\enddemo

We now show that the points in $A$ whose immediate stable 
manifolds 
are long are measure-theoretically dense 
in the circle $A$. This will imply that the basin of $A$ is 
measure-theoretically dense in the annulus and by a 
symmetric argument the 
basin of  
$B$ is also measure-theoretically dense and thus 
the two basins are intermingled.

\proclaim{Lemma 2.3} Let $R_{\rho}$  
denote the set $\{ {\bold x} \in A : \omega({\bold x}) = A, 
\,{\lambda}_{\perp}({\bold x}) < 0, \, \vert 
W^{\roman{is}}({\bold x}) \vert \geq 
\rho \} $. 
 Then, for 
any $0 < {\rho} <1$ we have $m_x(R_{\rho}) > 0$ and the 
points 
of $1$-dimensional Lebesgue density of $R_{\rho}$ are 
dense in $A$.
\endproclaim

Recall that ${\bold x}$ is a point of
Lebesgue density for a set $C$ if $\lim_{\varepsilon \to 0}
{{m({\Cal N}_{\epsilon}({\bold x}) \cap C)} \over
{m({\Cal N}_{\epsilon}({\bold x}))}} = 1$, where 
${\Cal N}_{\epsilon}({\bold x})$denotes an 
$\epsilon$-neighborhood 
of ${\bold x}$. Lebesgue almost every
point of a set of positive Lebesgue measure is a point of
Lebesgue density for that set.

\proclaim{Corollary} The basins of $A$ 
and $B$ are intermingled.
\endproclaim

\demo{Proof of Lemma {\rm2.3}} 
Let $R $ denote the set 
$\{ {\bold x} \in A : \omega({\bold x}) = A, 
\,{\lambda}_{\perp}({\bold x}) < 0 \}$.
By Lemma 2.1 
this is a set of full $1$-dimensional Lebesgue measure.
By Lemma 2.2 we can fix an $\varepsilon > 0$  so that  
$ m_x(R_{\varepsilon}) > 0 $.  

The immediate stable manifold of the
fixed point ${\bold p} = ({1\over 2},\,0)$ 
is the vertical line segment ${1 \over 2} \times [0,\,1)$ 
and is 
the saddle connection between the hyperbolic saddle fixed 
point
${\bold p}$ and 
the repelling fixed point ${\bold q}= ({1 \over 2},\,1)$.
Let $\theta(x,\,z)$ denote the branch of ${f}^{-1}$ 
which maps the annulus to the vertical strip $[{ 1\over 
3},\,
{2 \over 3}) \times [0,\,1].$  
By the vertical structure of $f$
we see that 
for each ${\bold x} \in {R_{\varepsilon}}$ the vertical
line segment   
$\theta^{k}(W^{\roman{is}}({\bold x})) $
converges to $W^{\roman{is}}({\bold p}) =\{{1 \over 2}\} 
\times 
[0,\,1)$ as $k \to \infty$.  
Let $n_{\delta}({\bold x})$ denote the smallest positive 
integer $n$ such
that  $\vert \theta^{n}(W^{is}({\bold x}))\vert \geq ( 
1-\delta)$. It is 
easy to see that for positive $ \delta$ the function 
$n_{\delta}({\bold x})$ is a bounded function on 
$R_{\varepsilon}$.
We denote the supremum of 
$n_{\delta}({\bold x)}$ on $R_{\varepsilon}$ by 
$n_{\delta}$.
 Since $\vert [D_{z}f_{z}](x,\,z) \vert \leq {33 \over 
32}$ for all 
$(x,\,z)$ in the annulus, for any 
$\theta^{n}(W^{\roman{is}}({\bold x})) $
whose length 
is at least $(1-\delta)$,  each of the $3^k$ evenly spaced 
(in the $x$-coordinate) images of 
$\theta^{n}(W^{\roman{is}}({\bold x})) $
under ${f^{-k}}$ is a vertical line segment of 
length at least $1-{({33 \over 32})}^{k}\delta$.  
     
So, for a fixed $0< {\rho}<1 $ we find a dense subset 
${\hat R}_{\rho}$ 
of $R_{\rho}$ with 
immediate stable manifolds of length greater than $\rho$ 
by taking  preimages of stable manifolds whose lengths 
are almost $1$. That is, for ${\rho}=1-\delta$ we have  
$${\hat R}_{\rho} = \bigcup_{ k\geq 0}
{f}^{-k}({\theta}^{n_{\ell(k)}}({R_{\varepsilon}})),$$
where $\ell(k) = \delta({32 \over 33})^k$.
Since $m_x(R_{\varepsilon})$ is positive, 
we see that 
$m_x\bigl({\theta}^{n_{\ell(k)}}({R_\varepsilon}) \bigr)$ 
is also positive.
The $3^k$ images under $f^{-k}$ of $R_\varepsilon$ 
also each have positive measure and 
are evenly spaced on $A$, and thus the set ${\hat 
R}_{\rho}$ has points of 
$1$-dimensional Lebesgue density which are dense in $A$. 
Since ${\hat R}_{\rho}$ is a subset of 
$R_{\rho}$, we see that $R_{\rho}$ is as desired.   
\enddemo

In order to show that the union of the basins of $A$ and 
$B$ has
full measure, we need to establish some control on the 
slope of iterated 
horizontal line segments.

\proclaim{Lemma 2.4} If $\gamma$ is a horizontal line 
segment of the form 
$\gamma = [n3^{-r},
(n+1)3^{-r}) \times \{ \rho \}$, then 
$f^r(\gamma) = {\gamma}_r$ is the graph of a function 
$Z(x)$ whose 
domain is $[0,\,1)$ and whose derivative is at most ${\pi 
\over 10}$ 
in absolute value.
\endproclaim

\demo{Proof of Lemma {\rm2.4}} The maximum value of 
$\vert [D_{z}{f}_{z}](x,\,z) \vert $
is ${33 \over 32}$, 
and the maximum value of $\vert [D_{x}{f}_{z}](x,\,z) 
\vert $ is
${\pi \over 16}$;
so the lower triangular form of the derivative matrix of 
(3) 
allows us to establish easily 
that slopes whose absolute value is less than ${\pi \over 
10}$  
continue to be slopes whose absolute value is less than 
${\pi \over 10}$ 
after iteration by $f$.  
\enddemo

\demo{Proof of Theorem $2$ restricted to the case of $f$}
By the corollary of Lemma 2.3 the basins of attraction of 
$A$ and
$B$ are intermingled.
In order to show that 
the complement of $\bigl( {\Cal B}(A) \cup {\Cal B}(B) 
\bigr)$ 
has $2$-dimensional Lebesgue measure zero, we will show 
that 
the complement of $\bigl( {\Cal B}(A) \cup {\Cal B}(B) 
\bigr)$  
has no points of horizontal $1$-dimensional Lebesgue 
density. 

Consider an arbitrary horizontal 
line segment $\gamma$ as in Lemma 2.4 and its $r$-th 
iterate ${\gamma}_{r}$.
We assume that $Z(0) \leq {1 \over 2}$.  
By Lemma 2.4
the derivative of the function $Z(x)$ is at most ${\pi 
\over 10}$, and thus the
function $Z(x)$ is bounded above by $\beta = {1 \over 2} + 
{\pi \over 10}$.
The immediate stable manifolds of the set $R_{\beta}$
each intersect ${\gamma}_{r}$ in exactly one point, so
$m_x\bigl(\gamma_r \cap {\Cal B}(A)\bigr) \geq 
m_x\bigl(R_{\beta}\bigr)$.  
The image under $f^{-r}$ in $\gamma$  of
$\gamma_r \cap {\Cal B}(A)$ is in the basin of attraction 
of $A$. Since $f^{r}$ is linear in the $x$-coordinate, 
$m_x\bigl(R_{\beta}\bigr)$ 
is a lower bound on the  $m_x$-measure of the fraction
of $\gamma$ which is in the basin of $A$. 
If $Z(0) > {1 \over 2}$, then we carry out a similar 
argument for $B$ instead 
of $A$. So the set of points which are in neither the 
basin of $A$ 
nor the basin of $B$ cannot have any points of horizontal 
$1$-dimensional Lebesgue density and 
by Fubini's theorem must 
therefore be a set whose $2$-dimensional measure is zero.
\qed
\enddemo

\heading 
3. The general case
\endheading 

We discuss the technical issues which must be dealt with 
to prove Theorem 2 
in the general case where ${\tilde f}$ is a small $C^{1+
\alpha}$ perturbation of $f$.  Full 
details can be found in \cite{6}. 
There are four main technical complications. 

First, the invariant measures on $A$ and $B$ may be 
perturbed. 
Fortunately, the perturbed measures remain absolutely 
continuous and have continuous densities 
close to 1. Second, if ${\tilde f}_x$ no longer 
depends only on the $x$-coordinate, vertical lines may 
no longer be sent to vertical lines. 
However, the results of Hirsch, Pugh, and Shub \cite{7} 
on normal hyperbolicity 
assure us that a kind of almost vertical structure will be 
preserved. 
Third, for a smooth curve $\eta$ transverse to the stable 
manifolds of $A$, 
the {\it holonomy} map $h_{\eta}$ taking a point ${\bold 
x} \in A$ 
to $\eta \cap W^{\roman{is}}({\bold x})$ is not 
necessarily even 
Lipschitz. Fortunately, 
P\'esin theory \cite{8--10} 
allows us sufficient control of the Radon-Nikodym 
derivative of 
the holonomy map. Fourth, distortion estimates are needed 
when 
pulling back measures from ${\gamma}_r$ to $\gamma$, since 
the map 
${\tilde f}$ need not be linear in the $x$-coordinate.

\heading Acknowledgment\endheading
\par We thank Charles Pugh, James C. Alexander, James A. 
Yorke, John Milnor,
Lai-Sang Young, and David Lavine for their helpful 
discussions.

\enddocument